\documentclass[a4paper,12pt]{article}
\usepackage{paralist}
\usepackage{amsmath,amssymb,amsfonts,amsthm}             

  \renewenvironment{thebibliography}[1]{%
    \begin{oldthebibliography}{#1}%
      \setlength{\parskip}{0.3ex}%
      \setlength{\itemsep}{0.0ex}%
  }%
  {%
    \end{oldthebibliography}%
  }

\usepackage[unicode, bookmarks, colorlinks, breaklinks,
	pdftitle={Thomas Carraro and Sven Wetterauer: XFEM Implementation},
	pdfauthor={Thomas Carraro and Sven Wetterauer: XFEM Implementation},
	pdfproducer={}
]{hyperref}  
\hypersetup{linkcolor=blue,citecolor=blue,filecolor=black,urlcolor=blue}

\usepackage[pdftex]{graphicx,color} 
\usepackage{subfig}
\usepackage{float}

\usepackage{lastpage}

\usepackage{numprint}

\usepackage{colortbl}
\usepackage{booktabs}
\definecolor{gray}{rgb}{0.7,0.7,0.7}

\definecolor{comment}{rgb}{.5,.5,.5}

\newtheorem{problem}{Problem}[section]

\begin{document}

\title{On the implementation of the eXtended Finite Element Method (XFEM) for interface problems}

\author{Thomas Carraro \thanks{thomas.carraro@iwr.uni-heidelberg.de} and Sven Wetterauer \thanks{sven.wetterauer@iwr.uni-heidelberg.de}\\Institute for Applied Mathematics, Heidelberg University}

\maketitle

\begin{abstract}
The eXtended Finite Element Method (XFEM) is used to solve interface problems with an unfitted mesh.
We present an implementation of the XFEM in the FEM-library deal.II.
The main parts of the implementation are (i) the appropriate quadrature rule; (ii) the shape functions for the extended part of the finite element formulation; (iii) the boundary and interface conditions.
We show how to handle the XFEM formulation providing a code that demonstrates the solution of two exemplary interface problems for a strong and a weak discontinuity respectively.
In the weak discontinuity case, the loss of conformity due to the blending effect and its remedy are discussed.
Furthermore, the optimal convergence of the presented unfitted method is numerically verified.
\end{abstract}

\section{Introduction}
\label{Introduction}
The eXtended Finite Element Method (XFEM) is a flexible numerical approach developed for {\bf general interface problems}.
Numerical methods to solve interface problems can be classified as \emph{fitted} or \emph{unfitted} methods. In the first case, the methods use a fitted mesh approach such that the interface is composed of element sides. The generation of a fitted mesh in case of complex interface geometry can be very time consuming. 
In many cases it cannot be done without handwork using a program for mesh generation.
In the unfitted case, the mesh is independent of the interface position and therefore unfitted methods are highly flexible. Since standard finite element methods perform poorly in the unfitted case, different alternative approaches have been introduced in the last years.

The XFEM is a partition of unity finite element method (PUFEM). The first formulation of the PUFEM has been derived in the work of Melenk and Babu{\u{s}}ka \cite{Melenk:1996}. The main features of this method can be summarized as follows:
\begin{compactitem}
\item a priori knowledge about the local behavior of the solution can be included in the formulation;
\item arbitrary regularity of the FE spaces can be constructed;
\item the approach can be understood as a meshless method;
\item it is a generalization of the $h$, $p$ and $hp$ version of the FEM.
\end{compactitem}
In particular two important aspects are essentially relevant for this method: local approximability and the capability to enforce inter-element continuity, i.e.\ conformity.
Among different PUFEM approaches, the generalized finite element method (GFEM) and the extended finite element method are the most versatile and the most used in many applications.
Their elaboration developed from the area of meshfree methods \cite{Belytschko:1996} and are based on the same principles: partition of unit and degree of freedom enrichment \cite{BelytschkoGracieVentura:2009}.
An overview of these methods can be found for example in \cite{Abdelaziz:2008,ANU:2003}.

The XFEM strategy to solve a problem with strong or weak discontinuities, i.\ e.\ discontinuities of the solution or of the fluxes respectively, is to extend the approximation space with discontinuous basis functions or basis functions with a kink, respectively.
Since the discontinuities are typically local features, the XFEM offers great flexibility by using a local modification of the standard FEM methods.
In fact, it avoids the use of complex meshing, which is substituted by a specific distribution of the additional degrees of freedom (DoFs).

The XFEM has broad use in different disciplines including fracture mechanics, large deformation, plasticity, multiphase flow, hydraulic fracturing and contact problems \cite{Khoei:2015}.
However, the first developments of the XFEM were done to simulate crack propagation \cite{Moes:1999}.
Further applications for the XFEM in material science comprise:
problems with complex geometries,
evolution of dislocations,
modeling of grain boundaries,
evolution of phase boundaries,
modeling of inclusions and homogenization problems.
In particular, the combination of the XFEM with a {\bf level set approach} \cite{Osher:2001} has been shown to be a very versatile tool to solve the above class of problems.
In the level set approach two strategies can be used to define the interfaces: (i) an analytical description of interfaces as the iso-zero of a function can be given or (ii) data from an image segmentation can be used to define interfaces locally or globally.

The goal of this work is to present all essential steps for the implementation of the XFEM. In addition, we make available a code (contact the authors to get a copy of it) that can be further extended for specific applications. The practical implementation is done in the open source FEM library deal.II \cite{dealII82}.
As application we consider an interface problem with a weak or a strong discontinuity. In the case of a strong discontinuity we consider only weak interface conditions of Robin type. In particular, we do not consider Dirichlet conditions on the interface. In this case the formulation of the problem has to be changed and possible formulations are based either on the Nitsche's method \cite{Nitsche:1971}, see for example \cite{Hansbo:2002,BeckerBurmanHansbo:2011}, or on a Lagrange multipliers method \cite{Moes:2006,Bechet:2009}. The extension of the code including the Dirichlet case is left for a further development.
In addition, we illustrate the problem of the \emph{blending effect}, i.e.\ of the loss of conformity in the elements adjacent to the interface, and show the implementation of a standard method to restore the full convergence behavior.
The focus of this work is on implementation aspects for the stationary XFEM.
We do not consider therefore moving interfaces. The development of a proper time stepping technique and an adequate quadrature rule goes beyond the scope of this paper.
We present examples in the two-dimensional case.
Specific aspects related to the extension to three-dimensional problems are also discussed.

This article is organized as follows. In section \ref{Interface problems} we introduce the general interface problem and the level set method. In section \ref{XFEM} the formulation of the XFEM for strong and weak discontinuities is depicted. Furthermore, the blending effect in case of weak discontinuities is discussed. In section \ref{A note on a priori estimation}, we briefly report some theoretical results on the convergence of the XFEM. In section \ref{Implementation} we describe in detail our implementation in deal.II. Specifically we discuss the XFEM quadrature rule and the application of boundary conditions to cut cells. In the final section \ref{Numerical examples} we show some numerical results on problems with weak and strong discontinuities including convergence tests.
\section{Interface problems}
\label{Interface problems}
\subsection{Problem setting}
\label{Problem setting}
Let $\Omega$ be a bounded domain in ${\mathbb R^2}$ with boundary $\partial \Omega$. The considered model problem is the stationary heat conduction. We construct an interface problem by taking a domain $\Omega$ divided in two sub-domains $\Omega_1$ and $\Omega_2$ by a line $\Gamma$, called {\em interface}.

\begin{figure}[h]
\centering
\scalebox{0.8}{
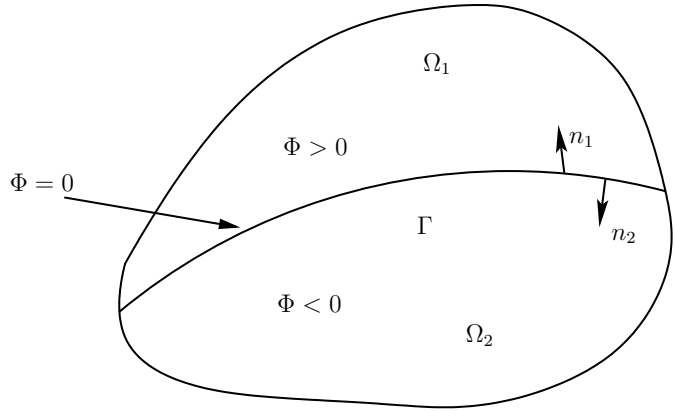
}
\caption{Level set function}\label{fig.domains}
\end{figure}

We consider three cases of interface problem with discontinuity on $\Gamma$. Case I: the solution has a discontinuity $g_s$; Case II: the solution has a weak discontinuity $g_w$; Case III: the solution has a discontinuity which depends on some functions $g_1$, $g_2$ as shown below.
The problem can be formulated as
\begin{problem}[Interface problem]
\label{interface problem}
Given the function $f$, the strong and weak discontinuity at $\Gamma$, $g_s$ and $g_w$ respectively, or the functions $g_1$, $g_2$, find the solution $u$ of the following system
\begin{subequations}
\begin{align}
-\nabla\cdot (\mu_i \nabla u_i) &= f &{\rm \,in\,} \Omega_i,\\
\label{dirichlet condition}
u_i &= g, &{\rm \,on\,} \partial\Omega_i \cap \Omega,\\
\text{Case I}:\quad\quad\quad\quad\quad\,
\label{strong discontinuity} [u] &= g_s &{\rm \,on\,} \Gamma,\\
\text{Case II}: \quad\quad\,\,
\label{weak discontinuity} [\mu \nabla u \cdot n] &= g_w &{\rm \,on\,} \Gamma,\\
[u] &= 0 &{\rm \,on\,} \Gamma,\\
\text{Case III}: \quad\quad
\label{weak discontinuity robin} \mu_i \nabla u_i \cdot n_i &= g_i(u_1, u_2) &{\rm \,on\,} \Gamma,
\end{align}
\end{subequations}
for $i=1,2$, with $[u] = u_1-u_2$ and $[\mu \nabla u \cdot n] = \mu_1 \nabla u_1 \cdot n_1 + \mu_2 \nabla u_2 \cdot n_2$, where $u_i$ is the restriction of $u$ to $\Omega_i$, $\mu_i$ is a constant assumed positive, and $n_i$ is the outward pointing normal to $\Omega_i$ at $\Gamma$, see Figure \ref{fig.domains}.
\end{problem}
We consider a general discretization with finite elements and use the notation with subscript $h$ to indicate discretized functions.
\begin{problem}[Discrete interface problem]
\label{discrete interface problem}
Given the same data as the continuous problem above, find the solution $u_h$ of the following discretized system
\begin{subequations}
\begin{align}
-\nabla\cdot (\mu_i \nabla u_{h,i}) &= f, &{\rm \,in\,} \Omega_i,\\
\label{dirichlet condition h}
u_{h,i} &= g, &{\rm \,on\,} \partial\Omega_i \cap \Omega,\\
\text{Case I}:\quad\quad\quad\quad\quad\,
\label{strong discontinuity h} [u_h] &= g_s &{\rm \,on\,} \Gamma,\\
\text{Case II}: \quad\quad\,\,
\label{weak discontinuity h} [\mu \nabla u_h \cdot n] &= g_w &{\rm \,on\,} \Gamma,\\
[u_h] &= 0 &{\rm \,on\,} \Gamma,\\
\text{Case III}: \quad\quad
\label{weak discontinuity robin h} \mu_i \nabla u_{h,i} \cdot n_i &= g_i(u_{h,1}, u_{h,2}) &{\rm \,on\,} \Gamma,
\end{align}
\end{subequations}
for $i=1,2$, with $[u_h] = u_{h,1}-u_{h,2}$ and $[\mu \nabla u_h \cdot n] = \mu_1 \nabla u_{h,1} \cdot n_1 + \mu_2 \nabla u_{h,2} \cdot n_2$, where $u_{h,i}$ is the restriction of $u_h$ to $\Omega_i$.
\end{problem}
Note that the same interface $\Gamma$ of the continuous problem is used also for the discretized problem if an exact quadrature formula can be employed.
An exact representation of $\Gamma$ used in the discretized problem can be obtained adopting the \emph{level set method} \cite{Osher:2001}.
\subsection{Level set method}
\label{Level set method}
The level set method is used to implicitly define the position of the interface $\Gamma$ independently of the underlying mesh used to discretize the interface problem.
The interface is defined by a scalar function $\phi: \Omega \rightarrow \mathbb R$ that is (uniquely) zero on $\Gamma$ and has different sign on different sub-domains:
\begin{align}
\phi&=0 &&\text{on }\Gamma,\notag\\
\phi&<0&&\text{in }\Omega_1,\\
\phi&>0&&\text{in }\Omega_2\notag.
\end{align}
Typically, the distance function (with sign) to the interface $\Gamma$ is used as level set function
\begin{align*}
\phi(x)=\pm \min_{y\in \Gamma}\Vert x-y\Vert.
\end{align*}
This is not the only choice, but it is often used because it can be exploited to calculate the normal vector at any point on $\Gamma$, since $\nabla \phi/|\nabla \phi|$ represents the normal vector if $\phi$ is the distance function to $\Gamma$.
Following a standard definition of XFEM we use the level set function to extend the space of finite elements with functions that incorporate the needed discontinuity.

\section{Extended finite elements}
\label{XFEM}
Generally, Galerkin finite elements are defined as the triple 
\begin{align}
({\cal T}, Q_p, \Sigma),
\end{align}
where $\cal T$ is the mesh, $Q_p$ is the space of test and trial functions and $\Sigma$ is a set of linear functionals that defines the degrees of freedom of the FEM formulation.
To build the XFEM in deal.II, we consider an extension of Lagrange finite elements. These are FE for which the degrees of freedom are the values of the test functions at the nodes of the mesh elements. Since our implementation is done in deal.II we consider only quadrilateral elements $K \in \cal T$.

The space $Q_p$ is the space of polynomial functions of degree at most $p$
\begin{align*}
Q_p=\{f(x)=\sum_{\alpha_1,\dots,\alpha_d\leq p}a_{\alpha}x_1^{\alpha_1}\dots x_d^{\alpha_d}  \}.
\end{align*}
defined on a unit cell $\hat K = (0,1)^d$. The test and trial functions on a real cell $K$ are obtained through a transformation $\sigma: \hat K \rightarrow K$ of a function from $Q_p$. We will use the notation $\varphi_{|_K} \in Q_p$ to indicate that the transformation of the function $\varphi_{|_K}$ onto the unit cell belongs to $Q_p$, i.e., $\sigma^{-1}(\varphi_K) \in Q_p$.

The scope of this work is to describe the implementation of the XFEM in deal.II, so we restrict our test cases to linear problems without loss of generality. We consider a general elliptic bilinear form and a linear functional:
\begin{align}
a: V\times V \rightarrow \mathbb R\\
f: V \rightarrow \mathbb R,
\end{align}
where $V$ is an appropriate Hilbert space.
The general weak formulation of our test cases is:
\begin{problem}
\label{continuous elliptic problem}
Find $u\in V$ such that
\begin{align}
a(u,\varphi) = f(\varphi), \quad \forall \varphi \in V.
\end{align}
\end{problem}
A typical choice for $V$ is the Hilbert space $H^1(\Omega)$ where $\Omega$ is the problem domain. 

The discrete approximation of problem \ref{continuous elliptic problem} using finite elements is
\begin{problem}
\label{discrete elliptic problem}
Find $u_h\in V_h$ such that
\begin{align}
a(u_h,\varphi_h) = f(\varphi_h), \quad \forall \varphi_h \in V_h,
\end{align}
\end{problem}
where $V_h$ is the finite dimensional space of $H^1$-conform functions
\begin{align*}
V_h&=\{\varphi\in V:\varphi\vert_K\in Q_p \text{ }\forall K\in \mathcal T \},
\end{align*}
The solution vector is a linear combination of the basis functions of $V_h$
\begin{align}
u_h(x)=\sum_{j=1}^n u_j N_j(x).
\end{align}
In the following we restrict our formulation to the space of bilinear functions $Q_1$, i.e.\ the shape functions $N_i$ are piecewise bilinear and globally continuous.

For interface problems of the type \ref{interface problem} there is the need to approximate a solution with a discontinuity. 
As illustrated above in the interface problem \ref{interface problem}, we consider three cases of boundary conditions on the interface leading to a weak discontinuity or a strong discontinuity.
In case of a weak discontinuity the standard $Q_1$ space can reach the best convergence rate (for a solution smooth enough) only if the mesh is fitted with the weak discontinuity.
In case of strong discontinuity, we consider convergence in a norm in the space $H^1(\Omega_1 \cup \Omega_2)$, since in the space $H^1(\Omega)$, where $\Omega=\Omega_1 \cup \Omega_2 \cup \Gamma$, the solution is not continuous and therefore it does not belong to $H^1$. Also in this case, the standard $Q_1$ space can achieve the best convergence rate (in $H^1(\Omega_1 \cup \Omega_2)$) only if the degrees of freedom lie on $\Gamma$.

If the interface cuts some elements $K$, a better approximation can be achieved by incorporating the discontinuity in the space in which we approximate the solution. 

\begin{center}
\begin{figure}[h!]
\begin{picture}(200,120)(-165,0)
\put(20,00){\line(1,0){120}}
\put(20,00){\line(0,1){120}}
\put(20,120){\line(1,0){120}}
\put(140,00){\line(0,1){120}}
\put(0,60){\line(1,0){160}}
\put(20,24){\line(1,0){120}}
\put(20,48){\line(1,0){120}}
\put(20,72){\line(1,0){120}}
\put(20,96){\line(1,0){120}}
\put(44,00){\line(0,1){120}}
\put(68,00){\line(0,1){120}}
\put(92,00){\line(0,1){120}}
\put(116,00){\line(0,1){120}}
\put(20,00){\circle*{5}}
\put(20,24){\circle*{5}}
\put(20,48){\circle*{5}}
\put(20,72){\circle*{5}}
\put(20,96){\circle*{5}}
\put(20,120){\circle*{5}}
\put(44,00){\circle*{5}}
\put(44,24){\circle*{5}}
\put(44,48){\circle*{5}}
\put(44,72){\circle*{5}}
\put(44,96){\circle*{5}}
\put(44,120){\circle*{5}}
\put(68,00){\circle*{5}}
\put(68,24){\circle*{5}}
\put(68,48){\circle*{5}}
\put(68,72){\circle*{5}}
\put(68,96){\circle*{5}}
\put(68,120){\circle*{5}}
\put(92,00){\circle*{5}}
\put(92,24){\circle*{5}}
\put(92,48){\circle*{5}}
\put(92,72){\circle*{5}}
\put(92,96){\circle*{5}}
\put(92,120){\circle*{5}}
\put(116,00){\circle*{5}}
\put(116,24){\circle*{5}}
\put(116,48){\circle*{5}}
\put(116,72){\circle*{5}}
\put(116,96){\circle*{5}}
\put(116,120){\circle*{5}}
\put(140,00){\circle*{5}}
\put(140,24){\circle*{5}}
\put(140,48){\circle*{5}}
\put(140,72){\circle*{5}}
\put(140,96){\circle*{5}}
\put(140,120){\circle*{5}}
\put(20,48){\circle{10}}
\put(20,72){\circle{10}}
\put(44,48){\circle{10}}
\put(44,72){\circle{10}}
\put(68,48){\circle{10}}
\put(68,72){\circle{10}}
\put(92,48){\circle{10}}
\put(92,72){\circle{10}}
\put(116,48){\circle{10}}
\put(116,72){\circle{10}}
\put(140,48){\circle{10}}
\put(140,72){\circle{10}}
\put(0,61){$\Gamma$}
\end{picture}
\caption{Single dots depicts normal degrees of freedom. Double dots depicts the extended degrees of freedom.}
\label{xfem_dofs}
\end{figure}
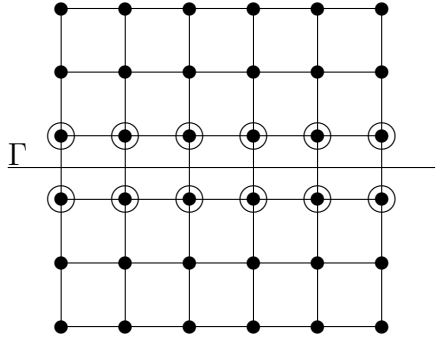

\end{center}
In the considered XFEM formulation we extend therefore the space $Q_1$ with some additional shape functions that represent the given discontinuity.
This is obtained by enriching the degrees of freedom of the elements cut by the interface. 
We will use the notation $I^\prime$ for the standard degrees of freedom and $I^*$ for the set of the extended degrees of freedom, respectively represented with single points and double points in Figure \ref{xfem_dofs}. 
The set of all degrees of freedom is denoted $I$.

In the next subsection we construct the XFEM shape functions.

\subsection{Strong and weak discontinuity}
In case of strong discontinuity, a typical function with a jump along $\Gamma$ is the sign function:
\begin{align*}
\text{sign}:\mathbb{R}&\rightarrow \{-1,0,1\}\\
\text{sign}(x)&=
\begin{cases}
1&\text{for }x>0\\
0&\text{for }x=0\\
-1&\text{for }x<0.
\end{cases}
\end{align*}
Since the jump is at the point $x=0$, the sign function applied to the level set function can be used to obtain a function with a jump along the interface.

In case of weak discontinuity, a function with a kink can be used, as for example the absolute value:
\begin{align*}
\text{abs}:\mathbb{R}&\rightarrow \mathbb{R}^+\\
\text{abs}(x)&=
\begin{cases}
x&\text{for }x>0\\
0&\text{for }x=0\\
-x&\text{for }x<0,
\end{cases}
\end{align*}
which applied to the level set function defines a weak discontinuity along the interface.

In the following we use the general notation
\begin{align*}
\psi:\Omega&\rightarrow \mathbb{R}\\
\psi(x)&=
\begin{cases}
\text{sign}(\phi(x))&\text{for strong discontinuity}\\
\text{abs}(\phi(x))&\text{for weak discontinuity},
\end{cases}
\end{align*}
$\psi$ is called {\em enrichment function}.

\begin{figure}[h!]
\begin{minipage}{0.5\textwidth}
\center
\includegraphics[scale=0.25]{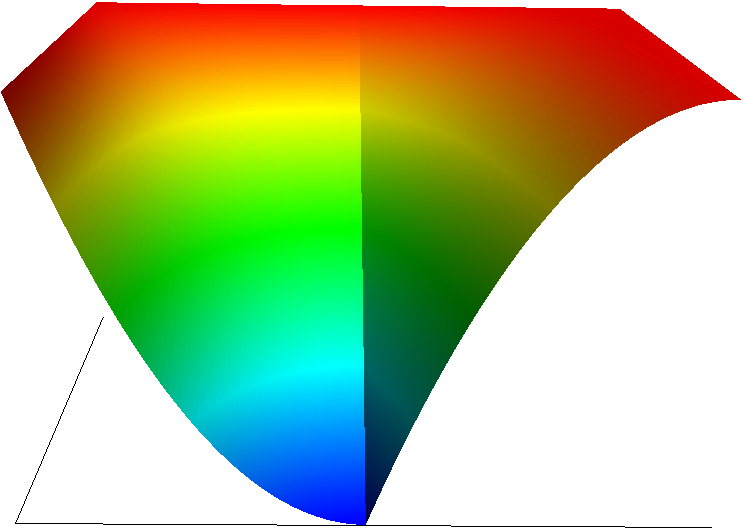}
\end{minipage}
\begin{minipage}{0.5\textwidth}
\center
\includegraphics[scale=0.25]{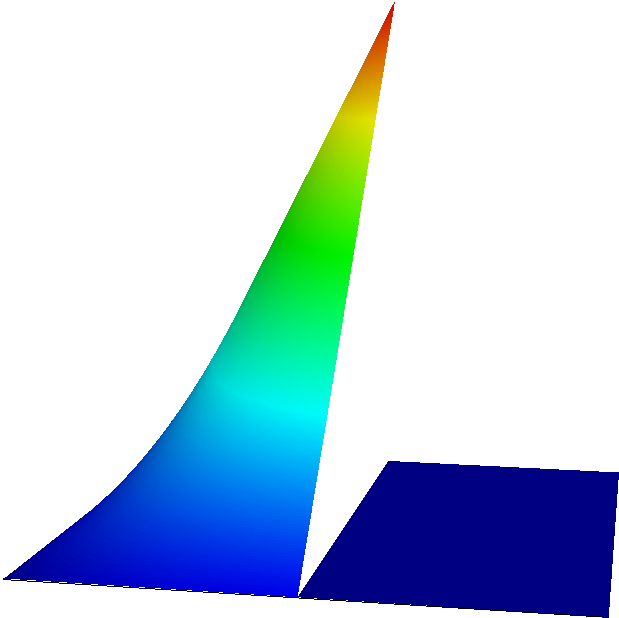}
\end{minipage}
\caption{XFEM basis functions for weak discontinuity (left) and strong discontinuity (right).}
\label{xfem basis}
\end{figure}

In Figure \ref{xfem_dofs} the extended degrees of freedom are depicted. These are additional Lagrangian degrees of freedom defined on a subset of existing mesh nodes.
To construct the XFEM, additional shape functions have to be defined. Following the above construction, we take functions that have a discontinuity along $\Gamma$
\begin{align}
M_i(x):=N_i(x)\psi(x).
\end{align}
We consider thus the following discrete spaces for test and trial functions
\begin{itemize}
\item strong discontinuity
\begin{align*}
V_h^s:=\{\varphi\in V:\varphi\vert_K\in Q_1\text{, }\varphi\vert_{K^\prime_i}\in Q_1, i=1,2\},
\end{align*}
where $K$ are standard cells and $K^\prime$ are the cells cut by the interface and $K^\prime_i$ are the two parts of the cut cell $K^\prime$ that have the interface as common edge.
\item weak discontinuity:
\begin{align*}
V_h^w:=\{\varphi\in V:\varphi\vert_K\in Q_1\text{, }\varphi\vert_{K'}\in Q_1\oplus \vert\phi\vert Q_1 \}.
\end{align*}
\end{itemize}
The discrete solution is now defined using the enriched basis
\begin{align*}
u_h(x)=\sum_{i\in I'} u_iN_i(x)+\sum_{j\in I^*} a_jM_j(x).
\end{align*}
From the practical point of view, it is desired that the discrete solution has the Kronecker delta property
\begin{equation}
u_{h}(x_{i})=u_{i}\quad i=1,\dots,n,
\label{kronecker} 
\end{equation}
where $x_i$ is the position of the $i^{th}$ degree of freedom.
To obtain this property, the extended basis functions are shifted, i.e. we use a different basis. The modified extended basis function of the node $i$ becomes
\begin{equation}
\label{shifted_basis}
M_{i}(x)=N_{i}(x)(\psi(x)-\psi(x_{i})).
\end{equation}
In Figure \ref{xfem basis} two extended basis functions are depicted, for the cases of weak and strong discontinuity.
The extended basis functions depend on the position of the interface. In case of weak discontinuity, the use of standard enriched basis functions can lead to the loss of conformity in the elements adjacent to the cut cells. This so called {\em blending effect} introduces a reduction of the convergence rate as it is shown later in the numerical experiments.

\subsection{Blending effect}
\begin{figure}[ht]
\centering
\includegraphics[scale=0.3]{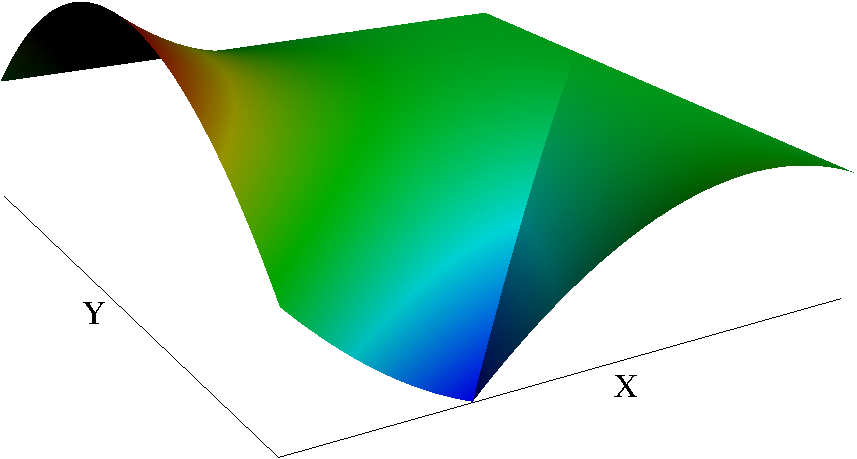}
\caption{Extended basis function with blending effect}
\label{blending}
\end{figure}

In this subsection we discuss the blending effect.
Let's consider the basis functions for the weak discontinuity
\begin{align}
\label{extended basis function}
M_i(x)=N_i(x)(\text{abs}(\phi(x))-\text{abs}(\phi(x_i))),
\end{align}
which is given by the product of a polynomial basis function $N_i$ and the level set function, which is in general a $C^2$ function. Considering the term $\text{abs}(\phi(x))-\text{abs}(\phi(x_i))$ along an edge $E$ of the cell where the function $N_i$ is not zero, it is
\begin{align}
\label{blending effect}
\text{abs}(\phi(x))-\text{abs}(\phi(x_i))
\begin{cases}
=\text{constant}&\text{for }E\parallel \Gamma\\
\neq \text{constant} &\text{for }E\nparallel\Gamma.
\end{cases}
\end{align}
Since the function $N_i$ along the edge $E$ is linear, the product with the enrichment function gives a function that is non linear. In Figure \ref{blending} an extended basis function is depicted that shows the blending effect.
The depicted function has a nonlinear behavior on the right side along the $x$ axis and on the left side along the $y$ axis. The nonlinearity along the $x$ axis is not a problem, since the edge is in common with a cut neighbor cell, which is also enriched in the same manner. On the other side, the edge along the $y$ axis causes problems, because the neighbor element is not cut by the interface and has only standard basis functions, which are linear on the common edge. We have thus linear behavior on one side and nonlinear on the other. Due to the discontinuity along this edge the enriched space $V_h^w$ is not $H^1$-conform anymore. Only in the case that the interface is parallel to the edges there is no blending effect, because the enrichment function is constant along such edges, see \eqref{blending effect}.

\begin{figure}[t!]
\center
\includegraphics[scale=0.3]{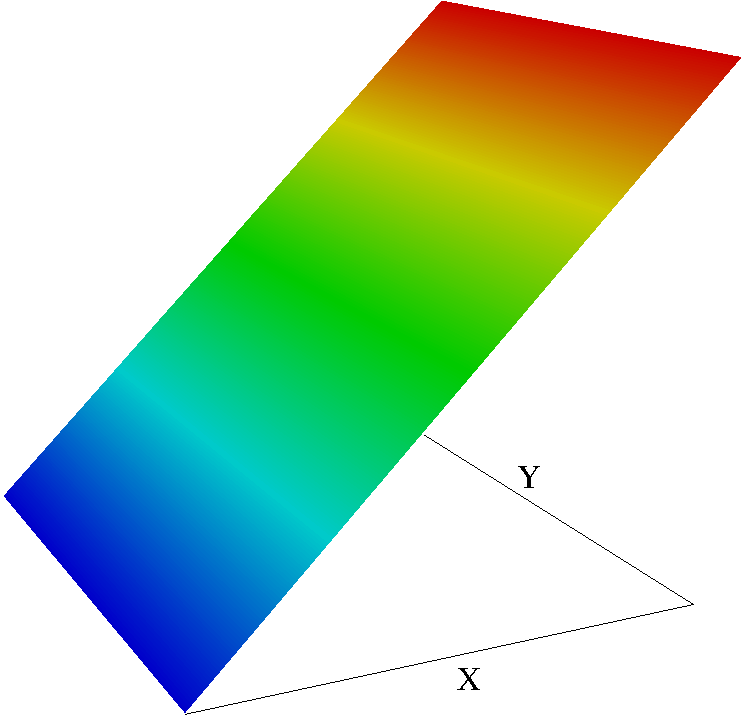}
\caption{Ramp function in a blending cell which neighbor cell on the right is a cut cell.}
\label{fig.ramp}
\end{figure}

\begin{figure}[ht]
\center
\includegraphics[scale=0.3]{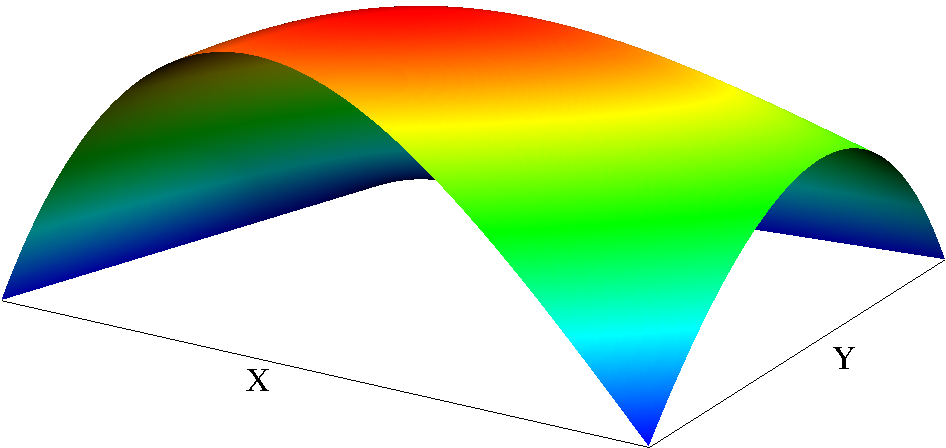}
\caption{Extended basis functions on blending cells.}
\label{blending_basis}
\end{figure}

There are different approaches to recover the $H^1$-conformity:
\begin{itemize}
\item use of higher order elements in the standard space $V_h$ \cite{Chessa:2003},
\item smoothing techniques as used in \cite{Sukumar:2000},
\item use of a corrected XFEM formulation adding a ramp function \cite{Fries:2008}. 
\end{itemize}
Note that the use of higher order elements works only if the extended functions on the cell edges have polynomial behavior.
In this work we use the third method using a correction with the ramp function:
\begin{align}
r(x):=\sum_{i\in I^\circ}N_i(x),
\label{ramp}
\end{align}
where $I^\circ$ depicts the set of standard degrees of freedom that lie on cut cells.
\\
The idea is to enrich not only the cells that are cut by the interface, but also the neighbor cells that are called {\em blending cells}.
In these cells the basis functions are modified so that they are nonlinear on the common edge with a cut cell and linear on the other edges recovering the global continuity.
\\
Figure \ref{fig.ramp} shows the ramp function on the neighbor of the cut cell depicted in Figure \ref{blending}.
The ramp function has the value 1 along the edge that creates the blending effect and is zero on the opposite edge. By multiplying the XFEM basis functions with the ramp function, new basis functions are defined that impose the continuity along the common edge with cut cells thus deleting the blending effect. Therefore, in the blending cells we use the extended functions:
\begin{align}
M_i(x):=N_i(x)(\text{abs}(\phi(x))-\text{abs}(\phi(x_i)))r(x).
\end{align}
which are depicted in Figure \ref{blending_basis}.
It can be observed that these additional functions behave nonlinearly along the common edge with a cut cell and a blending cell (in Figure \ref{blending_basis} depicted on the right side and on the bottom side respectively) and go to zero on common edges with normal cells (left and top side).
\section{A note on a priori error estimation}
\label{A note on a priori estimation}
In this section we briefly present some known results on a priori convergence estimates for the XFEM for interface problems. A result on optimal convergence rate of the XFEM for crack propagation can be found in \cite{Nicaise:2011}.
Let's consider the interface problem of the type III using the same notation as in Problem \ref{interface problem}
\begin{subequations}
\begin{align}
\label{pde_fuer_a_priori}
-\mu_i\Delta u&=f&&\text{in }\Omega_i\\
u&=0&&\text{on }\partial\Omega\\
\mu_1\partial_{n_1}u_1&=\alpha_{11}u_1+\alpha_{12}u_2&&\text{on }\Gamma\\
\mu_2\partial_{n_2}u_2&=\alpha_{21}u_1+\alpha_{22}u_2&&\text{on }\Gamma.
\end{align}
\end{subequations}
Note that continuity along the interface $\Gamma$ is not enforced.
For this problem we consider the weak formulation
\begin{problem}[Interface problem: Weak formulation]
Let the problem data be regular enough and let $\Omega_i$ be two domains with smooth boundaries so that the regularity $u\in H^2(\Omega_1\cup\Omega_2)$ is assured.
Find $u \in V$, such that for all $\varphi \in V$ it is
\begin{align}
a(u,\varphi)=(f,\varphi)_{\Omega_1\cup\Omega_2}
\end{align}
with 
\begin{align}
\label{robin-bilinear}
\begin{split}
a(u,\varphi)&=(\mu_1\nabla u,\nabla\varphi)_{\Omega_1}+(\mu_2\nabla u,\nabla\varphi)_{\Omega_2}\\
&\quad-(\alpha_{11}u_1,\varphi_1)_{\Gamma}-(\alpha_{12}u_2,\varphi_1)_{\Gamma}-(\alpha_{21}u_1,        \varphi_2)_{\Gamma}-(\alpha_{22}u_2,\varphi_2)_{\Gamma},
\end{split}
\end{align}
and $V=H_0^1(\Omega_1\cup\Omega_2;\partial\Omega)$, where we have used the notation $\partial\Omega:=\partial(\Omega_1\cup\Omega_2)\setminus \Gamma$.
\end{problem}
Let's consider the mesh $\mathcal T$ and the finite dimensional space $V_h\subset V$:
\begin{align}
\label{V_h}
V_h:=\{\varphi_h\in V:\varphi_h\vert_{K_i}\in Q_1, \quad \varphi_h\vert_{K_1} \equiv 0 \ \vee \ \varphi_h\vert_{K_2} \equiv 0, \quad \forall K\in \mathcal T, \quad i=1,2  \}
\end{align}
with $K_i:=K\cap \bar{\Omega}_i$. The basis functions $\varphi_h$ are the unfitted basis functions used in \cite{Hansbo:2002} to show the convergence results.
It can be shown that XFEM basis functions together with standard basis functions build a basis for the finite dimensional space $V_h$. In fact, as observed by Belytschko in \cite{comment_on_Hansbo} the unfitted basis functions in \cite{Hansbo:2002} are equivalent to the XFEM basis functions.
Therefore it is
\begin{align}
V_h=\text{span}\{N_i,\ M_j:\ i\in I^\prime,\ j\in I^* \},
\end{align}
where $I^\prime$ and $I^*$ are defined as in section \ref{XFEM}.
We consider the following XFEM approximation of the above problem

\begin{problem}[Interface problem: Discrete formulation]
With the same data as the above continuous problem, find $u_h\in V_h$, such that for all $\varphi_h\in V_h$ it is
\begin{align}
\label{discrete problem}
a(u_h,\varphi_h)=(f,\varphi_h)_\Omega,
\end{align}
with $\Omega = \Omega_1 \cup \Omega_2 \cup \Gamma$.
\end{problem}

Under these conditions following \cite{Hansbo:2002} it can be shown that for the finite element solution $u_h$ using the XFEM it is
\begin{align}
\Vert \nabla(u-u_h)\Vert_{\Omega_1\cup\Omega_2}\leq ch\Vert u\Vert_{H^2(\Omega_1\cup\Omega_2)},
\end{align}
and
\begin{align}
\Vert u-u_h\Vert_{\Omega_1\cup\Omega_2}\leq ch^2\Vert u \Vert_{H^2(\Omega_1\cup\Omega_2)}.
\end{align}
Note that to show these error estimations using the results in \cite{Hansbo:2002}, one has to perform a change of basis since Hansbo and Hansbo use different basis functions than the XFEM ones as pointed out above.

Therefore, full convergence behavior has to be expected in our numerical tests with strong discontinuity. Indeed as it will be shown later, also in case of weak discontinuity, using the ramp correction on blending cells, we observe the same convergence behavior. Nevertheless, the a priori convergence estimates for the weak discontinuity case cannot be derived in the same way following the work of Hansbo and Hansbo.
\section{Implementation in deal.II}
\label{Implementation}
To implement the XFEM in deal.II we consider the scalar interface problem \ref{interface problem} as a vectorial problem. The solution is therefore represented with two components. One component is the standard part of the solution and the other is the XFEM extension.
The standard part of the solution exists in all cells, whereas the extended part exists only in the cells that are cut by the interface, and their neighbors (blending cells) in case of weak discontinuity.
Due to an implementation constraint, both components of the vector-valued function must however exist in all cells. 
Therefore, since the degrees of freedom (in particular those belonging to the extended part) are distributed to all cells, we have to extend with the zero function the part of the extended solution over uncut cells.
By solving the system of equations it must be ensured that the zero extension of the solution remains also identically zero.

This is obtained in the implementation in deal.II with two main objects: the class {\em FeNothing} and the class {\em hp::DoFHandler}.
The object {\em FeNothing} is a finite element class that has zero degrees of freedom.
The object {\em hp::DoFHandler} allows to distribute different finite element types on different cells. Since we use the vector-valued finite element {\em FESystem} with two components, we can arbitrarily assign to each cell the two types of finite element either {\em FeNothing} or $Q_1$.
We assign thus a $Q_1$ finite element object to the first component (standard FE) of all cells, while we consider the following three cases for the second component of the {\em FESystem}:
\begin{itemize}
\item[(i)] for cells cut by the interface we use a $Q_1$ object to define the extended part of the space;
\item[(ii)] for the blending cells we use a $Q_1$ object to define the extended part of the space including the ramp functions.
\item[(iii)] for rest of the cells we use a {\em FeNothing} object (the extended part is set to zero);
\end{itemize}

The distribution of degrees of freedom with {\em hp::DoFHandler} is controlled by the value {\em active\_fe\_index} of the cell iterator.

\subsection{Quadrature formula}
\begin{figure}[ht]
\centering
\subfloat[\label{subfig.se1}]{
\includegraphics[scale=0.19]{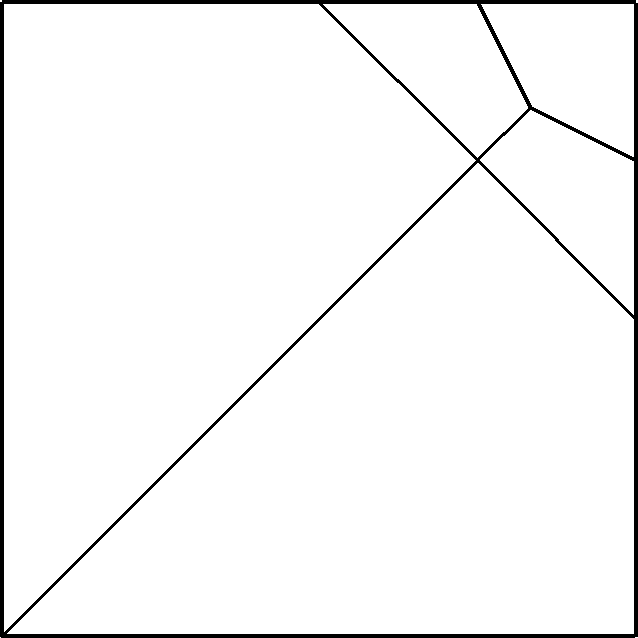}}
\subfloat[\label{subfig.se2}]{
\includegraphics[scale=0.19]{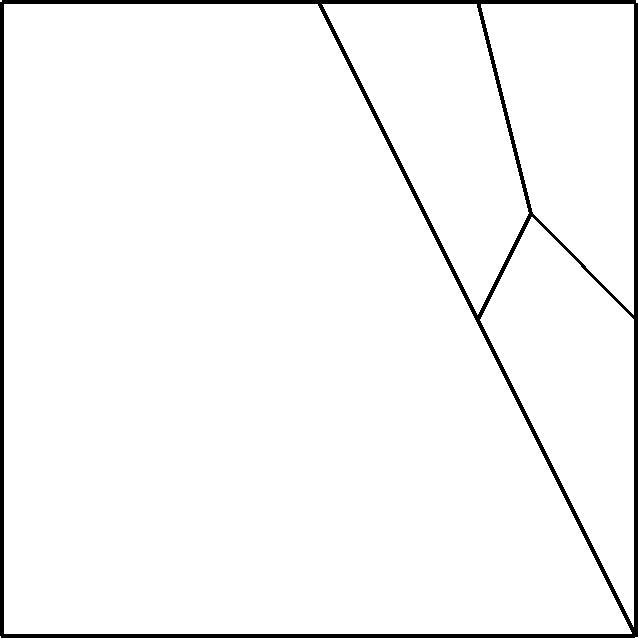}}
\subfloat[\label{subfig.se3}]{
\includegraphics[scale=0.19]{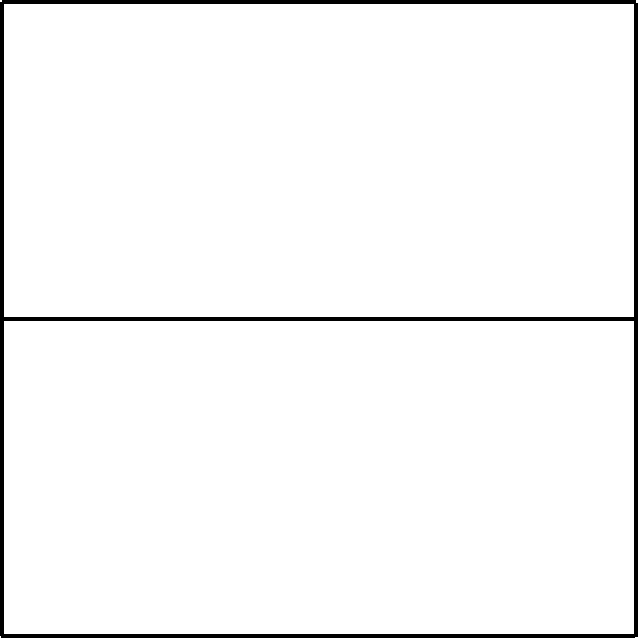}}
\subfloat[\label{subfig.se4}]{
\includegraphics[scale=0.19]{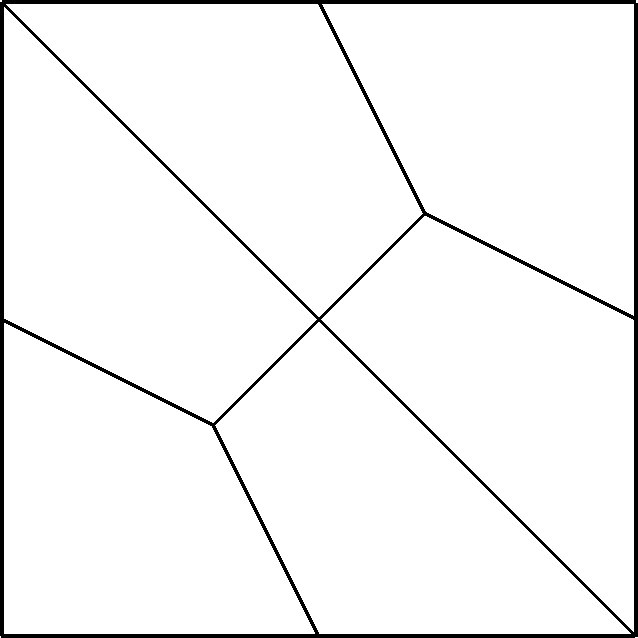}}
\caption{Subdivisions of the unit cell.}\label{2D subdivision}
\end{figure}

\begin{figure}[H]
\centering
\scalebox{0.5}{
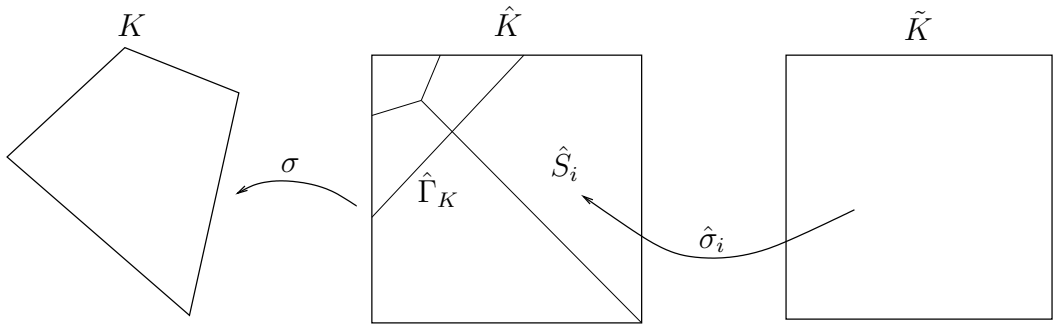
}
\caption{Transformations to construct the XFEM quadrature formula.}
\label{transformations}
\end{figure}

An essential part of the XFEM implementation is the quadrature formula. Generally, a quadrature formula such as the Gauss quadrature designed to integrate smooth functions would fail to integrate a function with a discontinuity.
A proper quadrature formula must take into account the position of the interface to integrate a function with a weak or strong discontinuity.
In our deal.II implementation this is done by subdividing the cut cells in sub-elements, on which a standard quadrature formula can be used.
Since we restrict our implementation to standard elements in deal.II, i.e. quadrilateral elements in 2D, the subdivision is done by quadrilateral sub-elements. In 2D there are only four types of subdivisions and the respective rotated variants, see Figure \ref{2D subdivision}, i.e. four of type (a), eight of type (b), two of type (c) and two of type (d).
Let's use the notation $K$ for the cell in real coordinates, $\hat K$ and $\tilde K$ for the unit cell, and  $\hat S_1, \dots, \hat S_n$ for the subcells of $\hat K$.
Since in deal.II the quadrature formula is defined for the reference unit cell, we have to construct a quadrature formula with points and weights to integrate a transformed function on the unit cell.
To this aim, the cut cell $K$ is transformed into the unit ``cut'' cell $\hat K$ and with the help of the transformed level set function it is subdivided according to the schemes in Figure \ref{2D subdivision}.
Each subcell is then transformed into the unit ``uncut'' cell $\tilde K$ in order to calculate the local quadrature points and weights through the standard tools in deal.II. 
Note that to make clear the use of the unit cell in different situations we use two notations for it, i.e. $\hat K$ and $\tilde K$.
Subsequently they are transformed back to the cell in ``real'' coordinates that in this case are the coordinates of the unit cell $\hat K$.
Different transformations are used to transform $K$ to $\hat K$ and each $\hat{S}_i$ into $\tilde K$ as depicted in the same figure. The transformation $\sigma$: $\hat K \rightarrow K$ is the standard transformation used in deal.II. 
Furthermore, to construct the appropriate XFEM quadrature formula we transform the subcell, see Figure \ref{transformations}, through the transformations:
\begin{align}
\hat \sigma_1,\dots,\hat \sigma_n:\tilde{K}\rightarrow \hat S_1,\dots, \hat S_n.
\end{align}
The XFEM quadrature formula is therefore derived by a standard quadrature formula.
For given quadrature points $x_i$ and weights $w_i$ of a standard formula, $i=1,\dots,m$, the XFEM quadrature on $\hat K$ is defined through the points
\begin{align}
y_{i,j}=\hat \sigma_j(x_i), \quad i=1,\dots,m,\quad j=1,\dots,n,
\end{align}
and weights
\begin{align}
w_{i,j}=w_i \, \text{det}\big(\nabla \hat \sigma_j(x_i)\big).
\end{align}
We define the degree of exactness of a quadrature formula as the maximal degree of polynomial functions that can be exactly integrated on an arbitrary domain. The optimal position of the quadrature points depends on the shape of the domain on which the integration is done.
The degree of exactness of the XFEM quadrature formula is at least of the same degree as the standard formula from which it is derived.
In addition, it allows by construction to integrate discontinuous functions along the interface.
From the implementation point of view this formula is highly flexible because it is built as a combination of standard formulas. The construction of an XFEM quadrature formula is simplified in our implementation since the subdivision of a cut cell is done by using the same type of cells, i.e.\ quadrilateral cells in our two-dimensional case.

\begin{figure}[H]
\centering
\subfloat[\label{subfig.xfem_quad1}]{
\scalebox{0.3}{
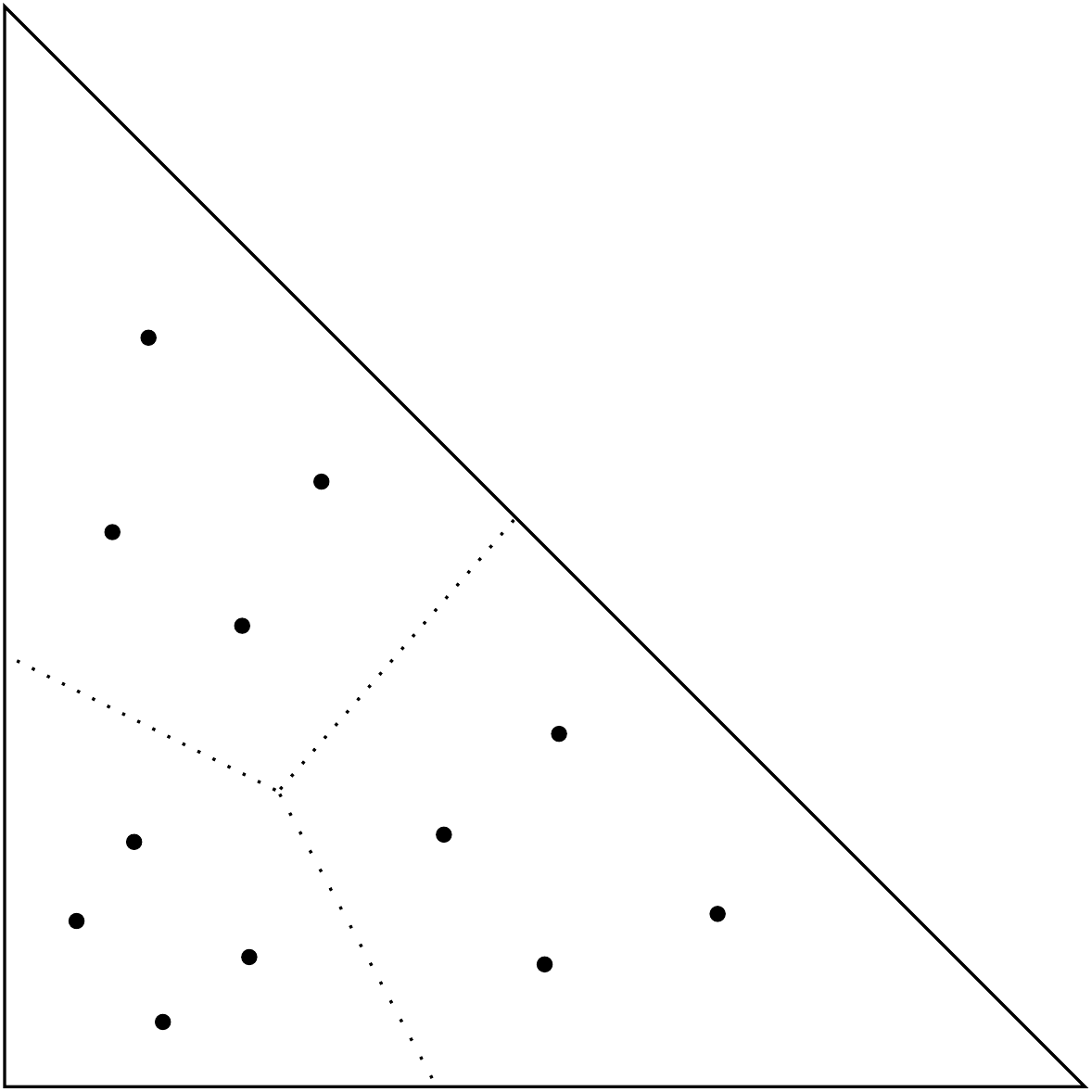
}
}
\hspace{2cm}
\centering
\subfloat[\label{subfig.xfem_quad2}]{
\scalebox{0.3}{
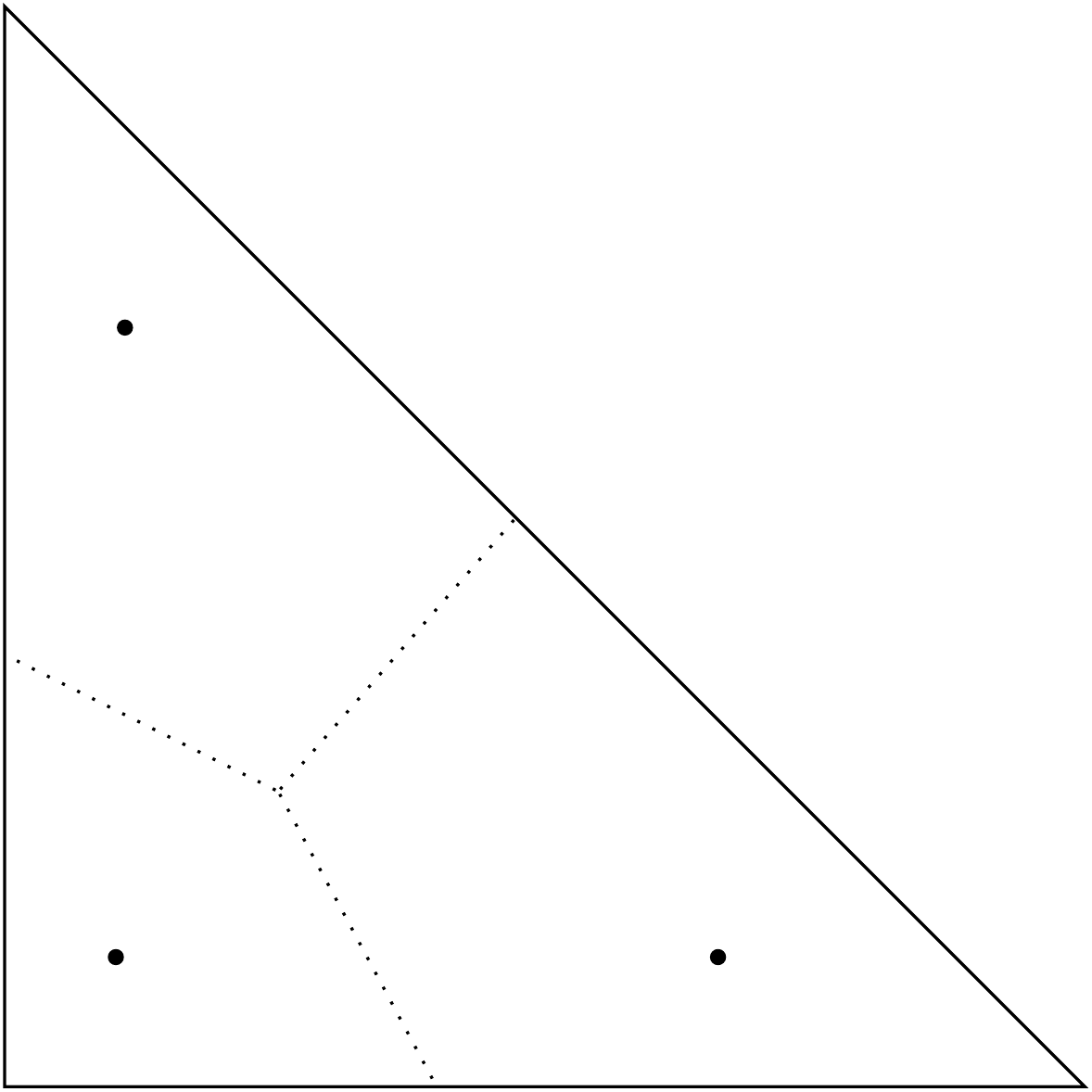
}
}
\caption{Sketch of the position in a triangle of the XFEM quadrature points (a) and of the symmetric Gauss formula (b).}
\label{XFEM quadrature triangle}
\end{figure}

\begin{figure}[H]
\centering
\subfloat[\label{subfig.3Dsubcell1}]{
\includegraphics[scale=0.18]{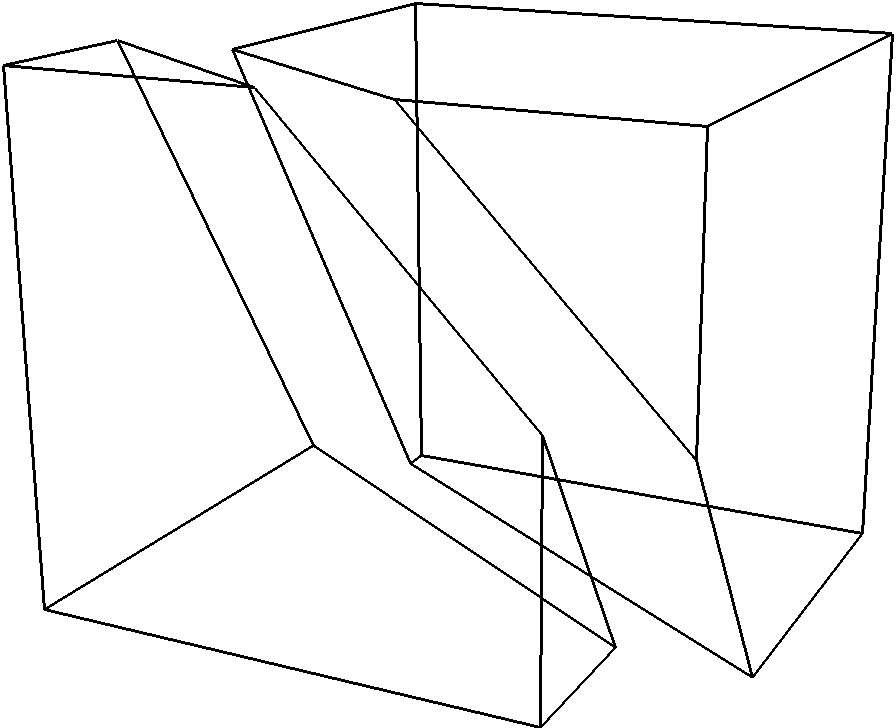}
}
\hspace{0.5cm}
\centering
\subfloat[\label{subfig.3Dsubcell2}]{
\includegraphics[scale=0.18]{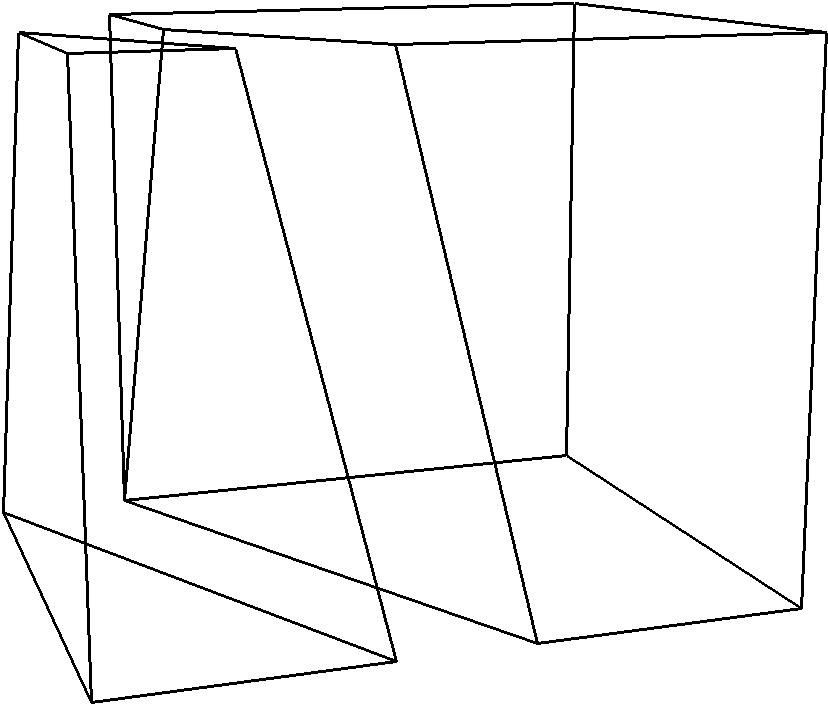}
}
\caption{Sketch of two partitions in the three dimensional case.}
\label{3D cuts}
\end{figure}

We would like to underline, that the XFEM formula is not optimal from the theoretical point of view, since it does not use the minimal number of points needed to integrate a given function over the specific subdivisions. In fact, the degree of exactness of the formula could be higher than the one inherited by the standard formula, but never less accurate. Therefore, the use of the XFEM formula can result in unnecessary higher costs to integrate a given function and therefore higher costs, for example, in assembling system matrices and vectors.
As an example, let's consider the integration of a quadratic function on the triangle resulting from the cut depicted in Figure \ref{2D subdivision} in the cases (a), (b) and (d).
The unit cut cell is divided in two parts, one triangle and one pentagon. For the pentagon part, quadrature rules would be necessary that are not typical on finite elements codes and therefore are not present as standard implementation. Therefore the division in two parts of the pentagon to obtain two quadrilaterals on which we can use standard formula is a desired feature. On the contrary on the triangle part one could use a more efficient formula, for example we could use a symmetric Gauss quadrature with 3 points, while the XFEM quadrature rule is build with 12 points as depicted in Figure \ref{XFEM quadrature triangle}. The code can be slightly more efficient changing the quadrature rule for the triangles obtained by the cut. We do not consider this modification for two reasons. On one side, we expect a great gain in computing time only in cases with a very high number of cut cells.
In a two dimensional numerical test not shown here we have observed about 13\% reduction of computing time using a three point quadrature formula for the triangle subdivisions instead of the XFEM one.
On the other, we want to produce a code that can be used in a dimension independent way.
In the two dimensional case the two parts of a cut cell can only be quadrilaterals, triangles or pentagons. In the three dimensional case there are 15 cases (and their respective symmetric configurations) and the partitions are more complex polyhedra as can be seen in the two depicted cases of Figure \ref{3D cuts}.
Therefore, in our implementation instead of considering a complex quadrature formula that can cope with all possible subdivisions, we apply subdivisions using only standard cells of deal.II in two and three dimensions, i.e.\ quadrilaterals and hexahedra respectively. The extension to the three dimensional case is theoretically straightforward. In practice the subdivision of hexahedral cut cells in hexahedral subcell is not a trivial task and it is left for a forthcoming work in which we will consider a comparison of different quadrature rules.

\subsection{Boundary and interface conditions}
This section is dedicated to the boundary and interface conditions. We consider Dirichlet and Neumann boundary conditions on the external boundary. Furthermore, we describe the implementation of Robin interface conditions on the interface $\Gamma$.
\subsection{Dirichlet and Neumann boundary conditions}
\label{sec.dirichlet}
In the following we consider the Dirichlet boundary condition \eqref{dirichlet condition} on the boundary of the domain. Nevertheless, the case with Neumann conditions can be treated in a similar way.

The Dirichlet boundary condition in deal.II is set by an appropriate modification of the system matrix and right hand side, and optionally performing one step of the Gaussian elimination process to recover original properties of the matrix as, e.\ g., the symmetry.

Owing to the Kronecker-delta property of the XFEM formulation, no special care has to be taken in case the Dirichlet condition is given by a continuous function. The degrees of freedom associated with the standard part of the FE are used to set the boundary values, while the degrees of freedom of the extension are set to zero at the boundary.
\begin{figure}[ht]
\center
\resizebox{0.9\textwidth}{!}{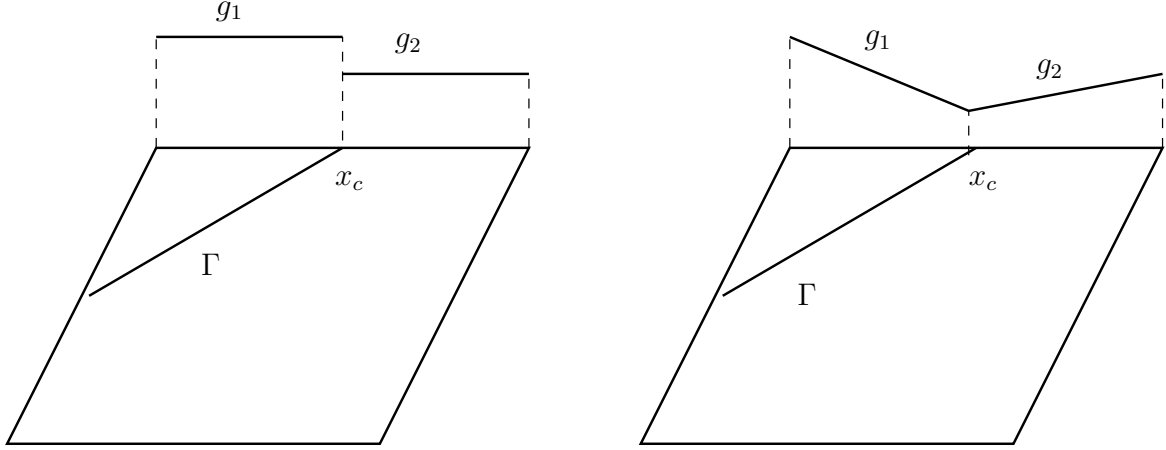}
\caption{Discontinuous boundary condition.}
\label{discontinuous BC}
\end{figure}
In case of discontinuous boundary condition, see Figure \ref{discontinuous BC} on the left side, the extended FE are used to approximate the discontinuity. On the edge at the boundary there are two extended degrees of freedom, whose basis function are discontinuous at one point. Due to the shift \eqref{extended basis function} used to construct the extended basis functions, these are nonzero only on one side of the edge. Therefore they can be used uncoupled to approximate the discontinuous Dirichlet value.

Let's consider the point $x_c$ of the discontinuity of $g$ on the edge with vertex $x_1$ and $x_2$, and the two limit values of the Dirichlet function $g_1(x_c)$ and $g_2(x_c)$. The two degrees of freedom $a_1$ and $a_2$ of the extended part at the boundary are uniquely defined by the following limit:
\begin{align}
\label{limit at x_c}
\lim_{x\to x_c,\ x\in\partial\Omega_i}u_h(x)=g_i(x_c),
\end{align}
which leads to
\begin{align}
a_1=\frac{g_2(x_c)-g_1(x_1)N_1(x_c)-g_2(x_2)N_2(x_c)}{2N_1(x_c)},\\
a_2=\frac{g_1(x_c)-g_1(x_1)N_1(x_c)-g_2(x_2)N_2(x_c)}{-2N_2(x_c)}.
\end{align}
In case of a weak discontinuity at the point $x_c$, the formulation needs a similar condition as in \eqref{limit at x_c} for the normal derivatives. The two extended degrees of freedom are therefore coupled and their value can be calculated solving a system of dimension $2\times2$.

\subsection{Robin interface conditions}
\label{sec.robin}
Following the notation of \eqref{weak discontinuity robin h} (strong discontinuity), to simplify the description of the interface conditions we assume that $g_1$ and $g_2$ are linear in both arguments.
For the variational formulation of the discrete interface problem \ref{discrete interface problem}, with interface conditions \eqref{weak discontinuity robin h}, we define the bilinear form
\begin{equation*}
  a(\varphi_h, \psi_h) = a(\varphi_h, \psi_h)_\Omega + a(\varphi_h, \psi_h)_\Gamma,
\end{equation*}
with
\begin{equation*}
  a(\varphi_h, \psi_h)_\Omega = (\mu_1\nabla \varphi_h,\nabla\psi_h)_{\Omega_1}+(\mu_2\nabla \varphi_h,\nabla\psi_h)_{\Omega_2}
\end{equation*}
and the boundary integrals
\begin{equation*}
  a(\varphi_h, \psi_h)_\Gamma = (g_1(\varphi_{h,1},\varphi_{h,2}),\psi_{h,1})_{\Gamma}+(g_2(\varphi_{h,1},\varphi_{h,2}),\psi_{h,2})_{\Gamma},
\end{equation*}
where the subscript 1 and 2 denotes the limit to the interface of the restriction on the subcells of the basis and test functions, e.g. for $\bar{x} \in \Gamma$ 
\begin{align}
  \label{limit shape function}
\varphi_{h,1}(\bar x)= \lim_{\Omega_1 \ni x \rightarrow \bar{x}}\varphi_h(x)
\end{align}
The bilinear form is used to build the system matrix and the scalar product has to be implemented considering all mixed products between standard and extended part of the finite element formulation.

In the case of continuous basis functions (standard FE part) $\varphi_{h,1}$ and $\varphi_{h,2}$ coincide.
On the contrary, in the extended part of the XFEM formulation the basis and test functions are discontinuous and one of the two functions in the scalar product vanishes due to the shift \eqref{shifted_basis}.
To determine the value of the limit \eqref{limit shape function} we use the {\em system\_to\_component\_index} in deal.II, which is a pair containing the component of the current DoF and the index of the shape function of the current DoF.
With the help of the level set function we can determine on which side of the interface the current DoF lies. This uniquely determines the zero part of the function.
\section{Numerical examples}
\label{Numerical examples}
In this section we consider the solution of a linear elliptic interface problem in case of strong and weak discontinuity. Furthermore, we show the effect of the blending cells
To simplify the notation in the following we use the symbol $\Omega$ instead of $\Omega_h$. The latter is the approximation of the boundary given by the adopted finite element formulation.

\subsection{Weak discontinuity}
\label{sec.weak}

\begin{figure}[t]
\center
\includegraphics[scale=0.5]{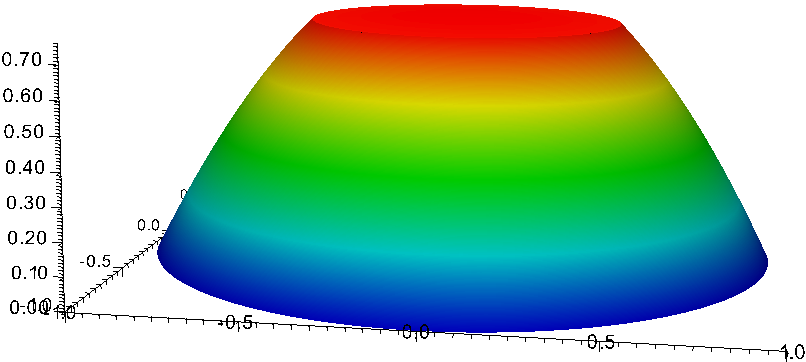}
\caption{Exact solution of problem \ref{pb.weak formulation of weak discontinuity}.}
\label{solution weak discontinuity}
\end{figure}

Let's consider the following domains $\Omega_1:=\{x\in \mathbb{R}^2:\Vert x\Vert_2 <0.5\}$, $\Omega_2:=\{x\in \mathbb{R}^2:0.5< \Vert x\Vert_2 <1\}$ and $\Omega:=\Omega_1\cup\Omega_2$.
with the interface
$\Gamma=\{ x\in \mathbb{R}^2\ : \Vert x \Vert _2 =0.5\}$.
We define the following problem
\begin{problem}[Weak discontinuity]
\label{pb.weak discontinuity}
Given $\mu_1=20$ and $\mu_2=1$, find the solution $u$ 
\begin{align}
\label{eq.weak_discontinuity}-\nabla\cdot(\mu_i\nabla u)&=1 && \text{in } \Omega_i,\\
\label{eq.weak_disc_dirichlet}u&=0 &&\text{on }\partial\Omega,\\
\label{eq.weak_disc_neumann}[u]&=0&&\text{on }\Gamma,\\
[\mu\partial_n u]&=0&&\text{on }\Gamma.
\end{align}
\end{problem}
\noindent The exact solution of problem \ref{pb.weak discontinuity} is the function $u_{ex}^w$
\begin{align}
u_{ex}^w:\Omega&\rightarrow \mathbb{R}\\
\notag
x& \mapsto 
\begin{cases}
\frac{1}{20}\cdot(-\frac{1}{4}\cdot\Vert x\Vert^2+\frac{61}{16}), & x\in\Omega_1,\\
\frac{1}{4}\cdot(1-\Vert x \Vert^2), & x\in\Omega_2.
\end{cases}
\end{align}
In Figure \ref{solution weak discontinuity} a XFEM approximation of this problem is shown.

\begin{problem}[Weak formulation of weak discontinuity]
\label{pb.weak formulation of weak discontinuity}
Find $u\in H^1_0(\Omega)$, so that 
\begin{align}
(\mu\nabla u,\nabla\varphi)_{\Omega}=(f,\varphi)_{\Omega} \quad \forall \varphi\in H^1_0(\Omega),
\label{formulierung_schwach}
\end{align}
with $\mu = \mu_1$ in $\Omega_1$ and $\mu = \mu_2$ in $\Omega_2$.
The conditions (\ref{eq.weak_disc_dirichlet}) -- (\ref{eq.weak_disc_neumann}) are naturally fulfilled by the weak formulation.
\end{problem}

\begin{figure}[ht]
\centering
\subfloat[\label{subfig.mesh1}]{
\includegraphics[scale=0.248]{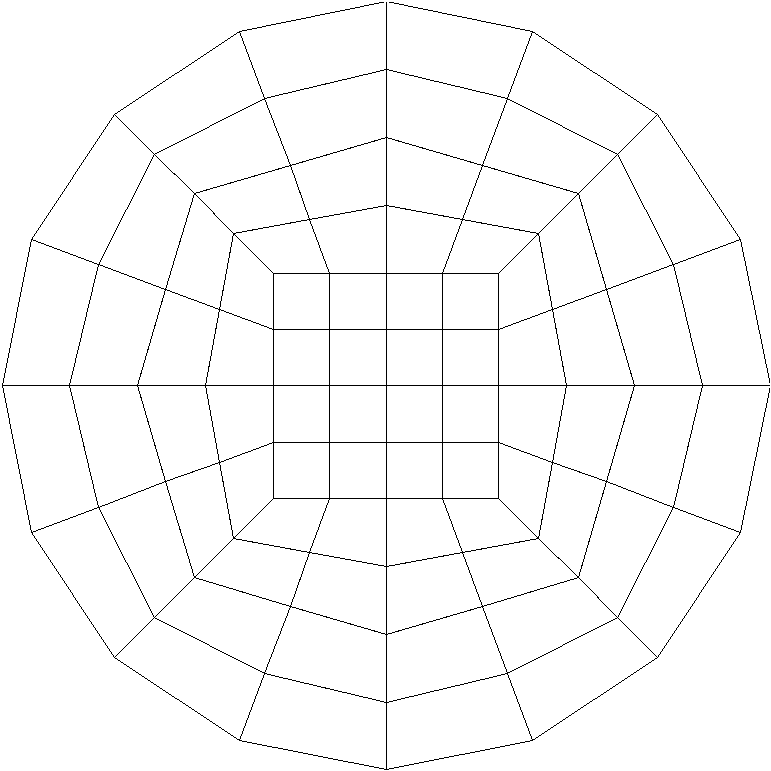}}
\subfloat[\label{subfig.mesh2}]{
\includegraphics[scale=0.3]{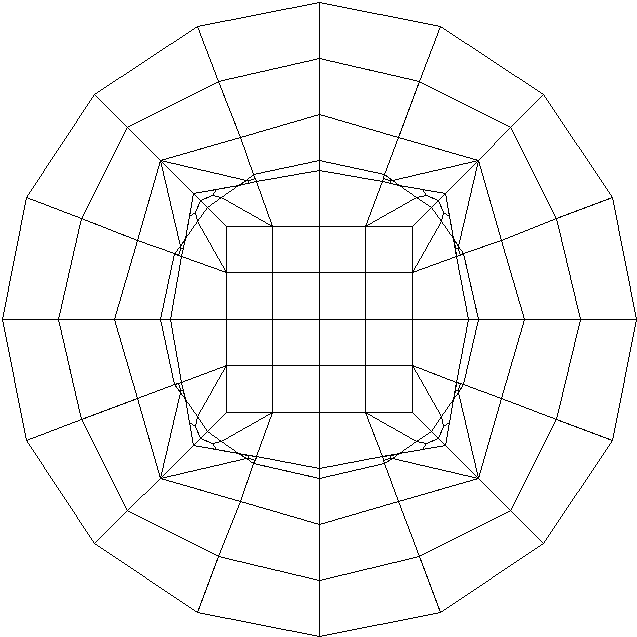}}
\subfloat[\label{subfig.mesh3}]{
\includegraphics[scale=0.3]{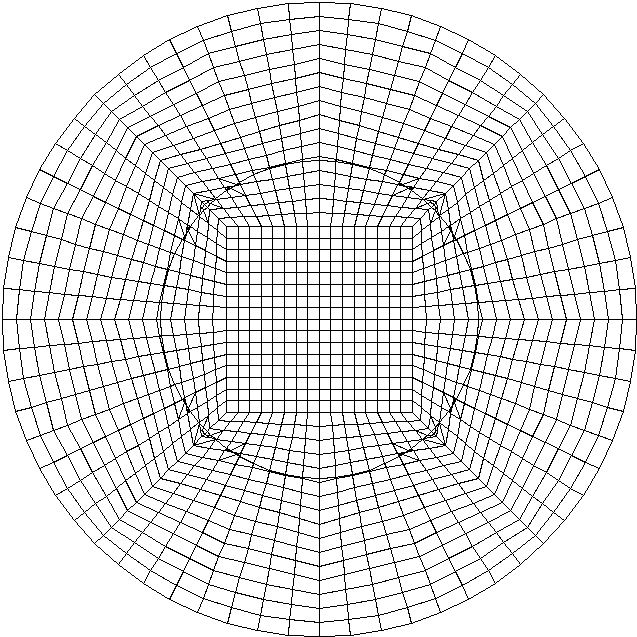}}
\caption{(a) Computing mesh at the coarsest level. (b) Coarsest mesh with visualized subcells; (c) Mesh at the second refinement level with visualized subcells.}\label{unfitted mesh}
\end{figure}

We consider in the following two approximations calculated using two different finite dimensional spaces to show the blending effect.
For a given mesh $\cal T$ let's consider the three subsets 
\begin{compactitem}
\item $\mathcal T_{cut}:=\{K\in\mathcal T:K\cap\Gamma\neq \emptyset \}$: the set of cells cut by the interface;
\item $\mathcal T_{bl}$: the set of blending cells, i.e. those cells that are neighbors of cut cells;
\item $\mathcal T_{std}:=\{K\in\mathcal T:K \not\subset \bigl({\mathcal T_{bl}} \cup {\mathcal T_{std}}\bigr)\}$: the set of standard cells.
\end{compactitem}
We can thereby define the $H_1$-conform space $V_h$ and the non-conform space $\tilde{V_h}$:
\begin{align*}
  \tilde{V_h}:=\{\varphi_h:&\varphi_h\vert_{K}\in Q_1\text{ for }K\in(\mathcal T_{std}\cup\mathcal T_{bl}),\\
&\varphi_h\vert_{K_i}\in Q_1\oplus|\phi| Q_1\text{ for }K\in\mathcal T_{cut},\\  
&\varphi_h \in C(\Omega\setminus\Gamma)  \}\\
V_h:=\{\varphi_h\in H^1_0(\Omega):&\varphi_h\vert_{K}\in Q_1\text{ for }K\in\mathcal T_{std},\\
&\varphi_h\vert_{K_i}\in Q_1\oplus|\phi| Q_1\text{ for }K\in\mathcal T_{cut},\\ 
&\varphi_h\vert_{K}\in Q_1\oplus r|\phi| Q_1\text{ for }K\in\mathcal T_{bl},\\ 
&\varphi_h \in C(\Omega)  \}.
\end{align*}
With this notation the discretized problem is given by:
\begin{problem}[Discrete formulation of weak discontinuity]
  With the data from Problem \ref{pb.weak formulation of weak discontinuity}, find $u_h\in V_h$, so that
\begin{align}
(\mu\nabla u_h,\nabla\varphi_h)_{\Omega}=(f,\varphi_h)_{\Omega} \quad \forall \varphi_h\in \hat V_h,
\end{align}
with either $\hat V_h=\tilde V_h$ or $\hat V_h=V_h$.
\label{discretized weak}
\end{problem}
We use an unfitted mesh, i.e.\ the interface $\Gamma$ intersects some cells. Therefore the underlying computing mesh is a shape regular mesh as shown in Figure \ref{unfitted mesh} (a).
Figures \ref{unfitted mesh} (b) and (c) show the subdivisions in subcells for two level of the mesh.
Remind that the subcells in the XFEM formulation are not finite element cells. They are only used to build the XFEM quadrature formula. In particular, even if the subcells are close to be degenerated this does not effect the quality of the mesh.

The numerical convergence results, shown in Table \ref{convergence weak} for the case with blending cells with ramp correction, show a full convergence rate. As expected for bilinear finite elements we observe a quadratic and linear convergence in the $L_2$ norm and energy norm respectively.
Table \ref{convergence weak no blending} shows the case without ramp correction for the blending cells. In this case, as expected, we observe a reduction of the convergence rate.
\begin{table}
\center
\begin{tabular}{r|r|r|r|r}
DoF & $\Vert u_h-u\Vert$ & Conv.rate & $\Vert \nabla(u_h-u)\Vert$ & Conv.rate\\
\hline
161		& 1.519e-02 & - &7.007e-02 & -\\
493		& 4.026e-03 & 1.92 &3.870e-02&	0.86\\
1621	& 9.980e-04 & 2.01 & 1.999e-02 & 0.95\\
5841 	& 2.533e-04 & 1.98 & 1.003e-02 & 1.00\\
21969 & 6.311e-05 & 2.00 & 5.026e-03 & 1.00\\
84945 & 1.575e-05 & 2.00 & 2.519e-03 & 1.00\\
\end{tabular}
\caption{Convergence rates of \ref{discretized weak} in $L_2$ and energy norm with conform space $V_h$.}
\label{convergence weak}
\end{table}
\begin{table}
\center
\begin{tabular}{r|r|r|r|r}
DoF & $\Vert u_h-u\Vert$ & Conv.rate & $\Vert \nabla(u_h-u)\Vert$ & Conv.rate\\
\hline
133 	& 1.155e-02 & - &1.002e-02 & -\\
417 	& 3.472e-03 & 1.73 &4.222e-02&1.25\\
1469 	& 8.671e-04 & 2.00 & 2.317e-02 & 0.87\\
5517 	& 2.214e-04 & 1.97 & 1.176e-02 & 0.98\\
21293 & 5.894e-05 & 1.91 & 6.559e-03 & 0.84\\
83565 & 2.007e-05 & 1.55 & 3.986e-03 & 0.72\\
\end{tabular}
\caption{Convergence rates of \ref{discretized weak} in $L_2$ and energy norm with non-conform space $\tilde{V_h}$.}
\label{convergence weak no blending}
\end{table}

\subsection{Strong discontinuity}
\label{sec.strong}

\begin{figure}[t]
\center
\includegraphics[scale=0.65]{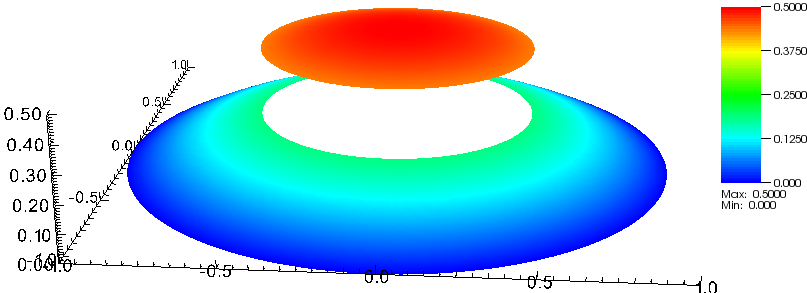}
\caption{Exact solution of problem \ref{pb.weak formulation of strong discontinuity}.}
\label{solution strong discontinuity}
\end{figure}

We consider the same domains as in the case of a weak discontinuity, $\Omega_1:=\{x\in \mathbb{R}^2:\Vert x\Vert_2 <0.5\}$, $\Omega_2:=\{x\in \mathbb{R}^2:0.5< \Vert x\Vert_2 <1\}$ and $\Omega:=\Omega_1\cup\Omega_2$.
On these domains we define the following problem:
\begin{problem}[Strong discontinuity]
\label{pb.strong discontinuity}
Find the solution $u$ 
\begin{align}
\label{eq.strong_discontinuity}-\Delta u&=1 && \text{in } \Omega_i,\\
\label{eq.strong_disc_dirichlet}u&=0 &&\text{on }\partial\Omega,\\
\label{eq.strong_disc_neumann}\nabla u_1\cdot n_1&=u_1-u_2&&\text{on }\Gamma,\\
\nabla u_2\cdot n_2&=-u_1+u_2&&\text{on }\Gamma.
\end{align}
\end{problem}
The exact solution of problem \ref{pb.strong discontinuity} is the function $u_{ex}^s$
\begin{align}
u_{ex}^s:\Omega_1\cup\Omega_2&\rightarrow \mathbb{R}\\
\notag
x& \mapsto
\begin{cases}
\frac{1}{4}\left ( 2-\Vert x \Vert_2 \right ), & x\in\Omega_1,\\
\frac{1}{4}\left( 1-\Vert x \Vert_2\right ), & x \in\Omega_2.
\end{cases}
\end{align}
In Fig \ref{solution strong discontinuity} a XFEM approximation of this problem is shown.

\begin{problem}[Weak formulation of strong discontinuity]
\label{pb.weak formulation of strong discontinuity}
Find $u\in H^1_0(\Omega_1\cup\Omega_2)$, so that
\begin{align}
(\nabla u,\nabla\varphi)_{\Omega}-(u_1-u_2,\varphi_1)_\Gamma-(-u_1+u_2,\varphi_2)_\Gamma=(f,\varphi)_{\Omega} \quad \forall \varphi\in H^1_0(\Omega_1\cup\Omega_2).
\label{formulierung_stark}
\end{align}
\end{problem}
Using the same notation as above the finite dimensional conform subspace of $H^1_0(\Omega_1\cup\Omega_2)$ is given by:
\begin{align*}
V_h:=\{\varphi_h\in H^1_0(\Omega_1\cup\Omega_2):&\varphi_h\vert_{K}\in Q_1\text{ for }K\in(\mathcal T_{std}\cup\mathcal T_{bl}),\\
&\varphi_h\vert_{K_i}\in Q_1\text{ for }K\in\mathcal T_{cut},\\
&\varphi_h\in C(\Omega\setminus\Gamma) \}.
\end{align*}
The discretized problem is then given by:
\begin{problem}[Discrete formulation of strong discontinuity]
Find $u_h\in\ V_h$, so that
\begin{align}
(\nabla u_h,\nabla\varphi_h)_\Omega-(u_{h,1}-u_{h,2},\varphi_{h,1})_{\Gamma}-(-u_{h,1}+u_{h,2},\varphi_{h,2})_{\Gamma}=(f,\varphi_h)_{\Omega} \quad \forall \varphi_h\in V_h.
\end{align}
\label{discretized_strong}
\end{problem}
The results of the numerical convergence analysis is given in Table \ref{convergence strong}.
Again, we observe the full convergence (quadratic in the $L_2$ norm and linear in the energy norm) in an unfitted mesh.
\begin{table}[h]
\begin{center}
\begin{tabular}{r|r|r|r|r}
DoF & $\Vert u_h-u\Vert$ & Conv.rate & $\Vert \nabla(u_h-u)\Vert$ & Conv.rate\\
\hline
133 	& 1.663e-02 & - &8.278e-02 & -\\
417 	& 4.503e-03 & 1.88 &4.231e-02&0.97\\
1469 	&1.045e-03 & 2.11 & 2.155e-02 & 0.97\\
5517 	& 2.875e-04 & 1.86 & 1.072e-02 & 1.01\\
21293 & 7.203e-05 & 2.00 & 5.364e-03 & 1.00\\
83565 & 1.803e-05 & 2.00 & 2.683e-03 & 1.00\\
\end{tabular}
\caption{Convergence rates of \ref{discretized_strong} in $L_2$ and energy norm.}
\label{convergence strong}
\end{center}
\end{table}

\section{Conclusions and possible extension of the program}

We have presented an implementation of the eXtended Finite Element Method (XFEM) in the FEM library deal.II.
The implementation is mainly based on the objects {\em hp::DoFHandler} and {\em FENothing}.

As part of this work, we make available a code that can be used to solve interface problems using the XFEM in two dimensions.
The main parts of the implementation are
\begin{compactitem}
  \item the XFEM quadrature formula;
  \item the assembling routine that uses the extended part of the finite element formulation;
  \item the visualization routine for cut cells.
\end{compactitem}
We have presented two prototypical examples to numerically approximate interface problems with strong and weak discontinuities respectively. The numerical results show the expected convergence rates. In case of a weak discontinuity also the blending effect (nonconformity) and the resulting loss of convergence rate are shown. Furthermore, we present the implementation of a known remedy to restore conformity.

We underly that a possible extension of this code is the implementation of a more efficient quadrature rule.
The implementation shown here can be straightforwardly extended to the three dimensional case by defining the necessary cell subdivisions.
In practice, this leads to a limited efficiency of the quadrature formula. Therefore, the focus of our next work is on an efficient quadrature rule for the three dimensional case.

An other extension of the program is the inclusion of extended shape functions for the approximation of crack propagation problems.
\begin{appendix}
  \section{Note on how to produce the numerical results}
  The source codes for two exemplary problems is included in the tar-file {\em xfem.tar}.
  Extracting this file two subfolders are produced, one for each problem. Note that the XFEM functions in the file {\em xfem\_functions.cc} are the same for both programs.
	\\
	To get the code running, one needs to adjust the {\em Makefile} in the according subfolder.
	In line 36 of the Makefile the correct path to the deal.II directory needs to be inserted.
	The used deal.II version needs to be at least version 8.0.
	The code is compiled with the command {\em make} and run with {\em make run}.
	\\
	Both programs will run performing several cycles of global mesh refinement.
	\subsection*{Strong discontinuity}
        The program in the subfolder {\em strong} produces the results shown in the Table \ref{convergence strong} of the Section \ref{sec.strong}. The following output is produced:
	\begin{verbatim}
Cycle 0:
   Number of active cells:       80
   Number of degrees of freedom: 133
   L2 error = 0.0166268
   energy error = 0.0827805
Cycle 1:
   Number of active cells:       320
   Number of degrees of freedom: 417
   L2 error = 0.00450276
   energy error = 0.0423088
Cycle 2:
   Number of active cells:       1280
   Number of degrees of freedom: 1469
   L2 error = 0.00104475
   energy error = 0.0215455
Cycle 3:
   Number of active cells:       5120
   Number of degrees of freedom: 5517
   L2 error = 0.000287494
   energy error = 0.0107175
Cycle 4:
   Number of active cells:       20480
   Number of degrees of freedom: 21293
   L2 error = 7.20349e-05
   energy error = 0.0053637
Cycle 5:
   Number of active cells:       81920
   Number of degrees of freedom: 83565
   L2 error = 1.80274e-05
   energy error = 0.00268312
       L2          Energy
1.663e-02    - 8.278e-02    -
4.503e-03 1.88 4.231e-02 0.97
1.045e-03 2.11 2.155e-02 0.97
2.875e-04 1.86 1.072e-02 1.01
7.203e-05 2.00 5.364e-03 1.00
1.803e-05 2.00 2.683e-03 1.00
	\end{verbatim}
	In the first part of the output information on the mesh and the errors in every refinement cycle is shown.
	In the last part of the output a table containing the $L_2$ and energy errors with the corresponding convergence rates is reported.
	\\
	Furthermore a file for the visualization of the solution is produced for every cycle in the vtk format.
	\subsection*{Weak discontinuity}
        The program in the subfolder {\em weak} can be used to solve an interface problem with weak discontinuities.
        There is an additional parameter file for the {\em weak} code, which contains the following parameters:
	\begin{verbatim}
	set Using XFEM        =true
	set blending          =true
	set Number of Cycles  =6
	set q_points          =3
	\end{verbatim}
	Four different parameters can be adjusted.
	In the first line the XFEM is activated.
        Setting this parameter to \verb!false! the program will use the standard FEM with unfitted interface.
	With the second parameter one can decide whether to use the ramp correction for the blending cells or not.
	Remind that this parameter has no effect, if the XFEM is not used.
	The third and fourth parameter set the number of refinement cycles and the number of quadrature points for each quadrature formula.
	\\
        This program produces the results shown in the Table \ref{convergence weak} of the Section \ref{sec.weak}. 
	Running it with the parameter set as above, the following output is produced:
	\begin{verbatim}
Cycle 0:
   Number of active cells:       80
   Number of degrees of freedom: 161
   L2 error = 0.01519
   energy error = 0.0700709
Cycle 1:
   Number of active cells:       320
   Number of degrees of freedom: 493
   L2 error = 0.0040257
   energy error = 0.0386981
Cycle 2:
   Number of active cells:       1280
   Number of degrees of freedom: 1621
   L2 error = 0.00099799
   energy error = 0.0199881
Cycle 3:
   Number of active cells:       5120
   Number of degrees of freedom: 5841
   L2 error = 0.000253274
   energy error = 0.0100287
Cycle 4:
   Number of active cells:       20480
   Number of degrees of freedom: 21969
   L2 error = 6.31119e-05
   energy error = 0.0050257
Cycle 5:
   Number of active cells:       81920
   Number of degrees of freedom: 84945
   L2 error = 1.57511e-05
   energy error = 0.00251584

       L2          Energy
1.519e-02    - 7.007e-02    -
4.026e-03 1.92 3.870e-02 0.86
9.980e-04 2.01 1.999e-02 0.95
2.533e-04 1.98 1.003e-02 1.00
6.311e-05 2.00 5.026e-03 1.00
1.575e-05 2.00 2.516e-03 1.00

	\end{verbatim}
	In the first part of the output information on the mesh and the errors in every cycle is shown.
	In the last part of the output a table containing the $L_2$ and energy errors with the corresponding convergence rates is reported.
        \\
	Furthermore a file for the visualization of the solution is produced for every cycle in the vtk format.
\end{appendix}

\section*{Acknowledgements}
We are thankful to Wolfgang Bangerth (Texas A\&M University) for helping with the initial implementation of the XFEM classes.
We acknowledge the work of Simon D\"orsam (Heidelberg University) for the visualization code for XFEM.
T.C.\ was supported by the Deutsche Forschungsgemeinschaft (DFG) through the project  ``Multiscale modeling and numerical simulations of Lithium ion battery electrodes using real microstructures (CA 633/2-1)''.
S.W.\ was supported by the Heidelberg Graduate School of Mathematical and Computational Methods for the Sciences.

\bibliographystyle{abbrv}
\bibliography{main}
\end{document}